 \documentclass[final,1p,times]{elsarticle}
\usepackage{amssymb}
\usepackage{amsmath}
\newtheorem{definition}{Definition}[section]
\newtheorem{theorem}{Theorem}[section]
\newtheorem{corollary}[theorem]{Corollary} 

\newtheorem{lemma}{Lemma}[section]

\numberwithin{equation}{section}
\journal{   }%Journal of Symbolic Computation
\begin{document}
	
	\begin{frontmatter}
		
		
		
		\title{Copositivity criteria of a class of fourth order 3-dimensional symmetric tensors} %% Article title
		
		
	\author[1]{Yisheng Song}
	
	\address[1]{School of Mathematical Sciences, Chongqing Normal University, Chongqing 401331 P.R. China.\\
		Email: yisheng.song@cqnu.edu.cn}
	\author[2]{Jinjie Liu}
	\address[2]{National Center for Applied Mathematics, Chongqing Normal University, Chongqing, P.R. China.\\
		Email: jinjie.liu@cqnu.edu.cn}
		
		%\affiliation{organization={School of Mathematical Sciences, Chongqing Normal University},%Department and Organization
			%addressline={}
			%city={Chongqing},
			%postcode={401331}, 
			%state={Chongqing},
			%country={China}}
		
		%% Abstract
		\begin{abstract}
			In this paper, we mainly dicuss  the non-negativity conditions for quartic homogeneous polynomials with 3 variables, which is the analytic conditions of copositivity of  a class of 4th order 3-dimensional symmetric tensors.  For a 4th order 3-dimensional symmetric tensor with its entries  $1$ or $-1$, an analytic  necessary and sufficient condition is given for its strict copositivity with the help of the properties of strictly semi-positive tensors. And by means of usual maxi-min theory,  a  necessary and sufficient condition is established for copositivity of such a  tensor also.  Applying these conclusions to a general 4th order 3-dimensional symmetric tensor,  the analytic conditions  are successfully obtained  for  verifying the (strict) copositivity, and these conditions can be very easily parsed and validated. Moreover, several  (strict) inequalities of ternary quartic homogeneous polynomial are established by means of these analytic conditions.
		\end{abstract}
		
	%	\begin{graphicalabstract}
		
	%	\begin{highlights}
		%	\item Research highlight 1
		%	\item Research highlight 2
		
		%% Keywords
		
	\end{frontmatter}
\footnotetext{The first author’s work was supported by the National Natural Science Foundation of P.R. China (Grant No.12171064), by The team project of innovation leading talent in chongqing (No.CQYC20210309536) and by the Foundation of Chongqing Normal univer- sity (20XLB009); The second author’s work was supported by the National Natural Science Foundation of P.R. China (Grant No.12001366), by the Science and Technology Research Program of Chongqing Municipal Education Commission(Grant No. KJZD-K202200506), by Chongqing Talent Program Lump-sum Project(No. cstc2022ycjh-bgzxm0040) and by the Foundation of Chongqing Normal University (22XLB005, ncamc2022-msxm02).}
% \linenumbers

\section{Introduction}

One of the most direct applications of $4$th order copositive tensors is to verify the vacuum stability of the Higgs scalar potential  model  \cite{K2012, K2016,CHQ2018,SL2021,S2022,S2023}.  In the graph theory, the $m$th order copositive tensors may be directly applied to estimate the bounds on the independent number of $m$-uniform hypergraph \cite{CHQ2018, NYZ2018, WCWYZ2023, CW2018}.
The concept of copositive tensors  was introduced by Qi \cite{Q2013} in 2013, which is usually applied to a symmetric tensor or, more precisely, to its associated Homogeneous polynomial of  degree $m$.  \begin{definition} Let $\mathcal{T}=(t_{i_1i_2\cdots i_m})$  be an
	$m$th order $n$ dimensional symmetric tensor.  $\mathcal{T}$ is called \begin{itemize}
		\item[(i)] {\bf positive semi-definite} (\cite{Q2005}) if $m$ is  an even number and  in the Euclidean space $ \mathbb{R}^n$, its associated Homogeneous polynomial $$\mathcal{T}x^m=\sum\limits_{i_1,i_2,\cdots,i_m=1}^nt_{i_1i_2\cdots
			i_m}x_{i_1}x_{i_2}\cdots
		x_{i_m}\geq0;$$ 
		\item[(ii)] {\bf positive definite} (\cite{Q2005}) if $m$ is  an even number and  $\mathcal{T}x^m>0$ for all $x\in \mathbb{R}^n\setminus\{0\}$;
		\item[(iii)] {\bf copositive} (\cite{Q2013}) if  $\mathcal{T}x^m\geq0$ on the nonnegative orthant  $\mathbb{R}^n_+$;
		\item[(iv)] {\bf strictly copositive} (\cite{Q2013}) if  $\mathcal{T}x^m>0$ for all $x\in \mathbb{R}^n_+\setminus\{0\} $.
\end{itemize} \end{definition}
Clearly, the positive semi-definite tensors must be copositive, and a copositive tensor coincides with a copositive matrix if $m=2$.   The concept of copositive matrices was introduced by   Motzkin \cite{M1952} in 1952.  Baston \cite{B1969} gave an analytic  way of judging copositivity of a $n\times n$ matrix in 1969.

\begin{theorem}(Baston \cite{B1969})\label{thm:11}
	Let $M=(m_{ij})$ be a symmetric  matrix with $|m_{ij}|=m_{ii}=1$ for all $i,j\in\{1,2,\cdots,n\}$. Then  the matrix $M$ is copositive if and only if there is no triple $(r,s,t)$ such that $$m_{rs}=m_{rt}=m_{st}=-1.$$
\end{theorem}
Simpson-Spector \cite{SS1983} and  Hadeler \cite{H1983} and Nadler \cite{N1992} and Chang-Sederberg \cite{CS1994} and Andersson-Chang-Elfving \cite{ACE}  respectively showed the (strict) copositive conditions  of $2\times2$ and $3\times3$ matrices  using different methods of argumentation.

\begin{theorem}\label{thm:12} Let $M=(m_{ij})$ be a symmetric  matrix. Then 
	a $2\times 2$ matrix $M$ is (strictly) copositive if and only if
	$$
	m_{11}\geq0\ (>0), m_{22}\geq0\ (>0), m_{12}+\sqrt{m_{11}m_{22}}\geq0\ (>0);
	$$
	a $3\times 3$ matrix $M$ is (strictly) copositive  if and only if for all $i\in\{1,2,3\}$,
	$$\begin{aligned}
		m_{ii}\geq 0\ (>0),
		\alpha=m_{12}+\sqrt{m_{11}m_{22}}\geq 0\ (>0),\\ \beta=m_{13}+\sqrt{m_{11}m_{33}}\geq 0\ (>0), \gamma=m_{23}+\sqrt{m_{33}m_{22}}\geq 0\ (>0),\\
		m_{12}\sqrt{m_{33}}+m_{13}\sqrt{m_{22}}+m_{23}\sqrt{m_{11}}+\sqrt{m_{11}m_{22}m_{33}}+\sqrt{2\alpha\beta\gamma}\geq 0\ (>0).
	\end{aligned}$$
\end{theorem}
Schmidt-He$\beta$ \cite{SH1988} provided the nonnegative conditions of a cubic and univariate polynomial  with real coefficients in non-negative real number $\mathbb{R}_+$, and Qi-Song-Zhang \cite{QYZ} recently gave the positivity of such a cubic polynomial,   which actually gave a the (strict) copositivity of 3rd order 2-dimensional symmetric tensor (see Liu-Song \cite{LS2021} for more details also).  Qi-Song-Zhang \cite{QYZ} also gave  the nonnegativity and positivity of a quartic and univariate polynomial   in $\mathbb{R}$, which means  the positive (semi-)definitiveness of 4th order 2-dimensional  tensor. Ulrich-Watson \cite{UW} and Qi-Song-Zhang \cite{QYZ} presented the analytic  conditions of the nonnegativity of a quartic and univariate polynomial   in $\mathbb{R}_+$.   This actually yielded   the copositivity of 4th order 2-dimensional symmetric tensor \cite{SL2021}. 

\begin{theorem} \label{thm:13} A 3rd order $2$-dimensional symmetric tensor $\mathcal{T}=(t_{ijk})$ is (strictly) copositive if and only if $t_{111}\geq0\ (>0)$,\ $t_{222}\geq0\ (>0)$, either $t_{112}\geq0,\ t_{122}\geq0$ or
	$$4t_{111}t_{122}^3+4t_{112}^3t_{222}+t_{111}^2t_{222}^2-6t_{111}t_{112}t_{122}t_{222}-3t_{112}^2t_{122}^2 \geq0\ (>0).$$
	A 4th order 2-dimensional symmetric tensor $\mathcal{T}=(t_{ijkl})$ with $t_{1111}>0$ and $t_{2222}>0$ is copositive if and only if
	$$\begin{cases}
		\Delta\leq0, t_{1222}\sqrt{t_{1111}}+t_{1112}\sqrt{t_{2222}}>0;  \\
		t_{1222}\geq0, t_{1112}\geq0, 3t_{1122}+\sqrt{t_{1111}t_{2222}}\geq0;  \\
		\Delta\geq0,\\
		|t_{1112}\sqrt{t_{2222}}-t_{1222}\sqrt{t_{1111}}|\leq\sqrt{6t_{1111}t_{1122}t_{2222}+2t_{1111}t_{2222}\sqrt{t_{1111}t_{2222}}}\\ 		
		(i) \ -\sqrt{t_{1111}t_{2222}}\leq3t_{1122}\leq3\sqrt{t_{1111}t_{2222}};\\ 		
		(ii) \ t_{1122}>\sqrt{t_{1111}t_{2222}} \mbox{ and }\\ 		
		t_{1112}\sqrt{t_{2222}}+t_{1222}\sqrt{t_{1111}}\geq-\sqrt{6t_{1111}t_{1122}t_{2222}-2t_{1111}t_{2222}\sqrt{t_{1111}t_{2222}}},
	\end{cases}$$
	where $
	\Delta =4\times12^{3}(t_{1111}t_{2222}-4t_{1112}t_{1222}+3t_{1122}^{2})^{3}
	-72^{2}\times6^{2}(t_{1111}t_{1122}t_{2222}+2t_{1112}t_{1122}t_{1222}-t_{1122}^{3}-t_{1112}^{2}t_{2222}-t_{1111}t_{1222}^{2})^{2}.
	$
\end{theorem}
For  a special 4th order 3-dimensional tensor given by the particle physical model, Qi-Song-Zhang \cite{QSZ2020} presented a necessary and sufficient condition of copositivity, and Song-Li \cite{SL2021} provided an analytic condition of its copositivity. However, an analytic necessary and sufficient condition  has not been found for  the copositivity of a general 3-dimensional higher order tensor ($m>2$). Even for a general 2-dimensional higher order tensor ($m>3$), people still has not found the analytic  conditions of strict copositivity untill now.  

For  checking the copositivity of symmetric tensors, various numerical algorithms have been  employed. For example, Chen-Huang-Qi \cite{CHQ2017,CHQ2018} gave the detection algorithms based on simplicial partition;   Li-Zhang-Huang-Qi \cite{LZHQ} used an SDP relaxation algorithm;  Nie-Yang-Zhang \cite{NYZ2018} devised a complete semi-definite algorithms.   Taking advantage of the properties of a tensor itself, the copositivity  can be described qualitatively.  Song-Qi \cite{SQ2015} showed a necessary and sufficient condition of copositivity by means of the principal sub-tensors; Song-Qi \cite{SQ2016} applied the sign of its Pareto H-eigenvalue (Z-eigenvalue) to test the copositivity. Song-Qi \cite{S-Q2016} proved the equivalence of (strict) copositivity and (strict) semi-positivity of a symmetric tensor. 
\begin{theorem}(Song-Qi \cite{S-Q2016})\label{thm:14}
	Let $\mathcal{T}=(t_{i_1i_2\cdots i_m})$ be a symmetric  tensor. Then   $\mathcal{T}$ is (strictly) copositive if and only if $\mathcal{T}$ is (strictly) semi-positive, i.e.,  for each $x=(x_1,x_2,\cdots,x_n)^\top\in \mathbb{R}^n_+\setminus\{0\}$, there exists $k\in\{1,2,\cdots, n\}$ such that $$x_k>0\mbox{ and }(\mathcal{T}x^{m-1})_k=\sum_{i_2\cdots i_m=1}^nt_{ki_2\cdots i_m}x_{i_2}\cdots x_{i_m}\geq0\ (>0).$$
\end{theorem}

For more  details of the copositivity of a higher order tensor and a matrix, see \cite{WCWYZ2023,  LV2022,HN1963, D1962,B1966,B1967,HS2010,  QCC2018,QL2017,DC2020}. Since there are  strong connection between the semi-positivity of a tensor and the tensor complementarity problems \cite{SQ2015,HQ2019,QH2019,H-Q2019}, so the copositivity of a symmetric tensor may be verified by solving the  corresponding tensor complementarity problems. For more details about the tensor complementarity problems and its applications, see Refs. \cite{MHW2023,SP2023,LT2023,BP,BHW,CQW2019,CQW2016,CW,CQS,G,HQ,LQN,SY,SQ2017,S-Q2015,SM2018,WCW,WHB,WHH}.

Motivation to checking  the copositivity of a higher order tensor,  we mainly dicuss analytic  necessary and sufficient conditions of copositivity of a class of 4th order 3-dimensional symmetric tensors in this paper. With the help of Theorem \ref{thm:14}, we  first promote Theorem \ref{thm:11} to ones of 4th order 3-dimensional symmetric tensors,  which gives  an analytic  necessary and sufficient condition  of strict copositivity of such a class of tensors (Theorem \ref{thm:34}). Secondly,  we present an analytic  necessary and sufficient condition  of copositivity of such a class of tensors (Theorem \ref{thm:38}).  Finally, applying Theorems \ref{thm:34} and \ref{thm:38},  the analytic sufficient conditions (Theorems \ref{thm:39}, \ref{thm:310} and \ref{thm:311})  are successfully proved for (strict) copositivity of a general 4th order 3-dimensional symmetric tensor. Furthermore,   several  (strict) inequalities of ternary quartic homogeneous polynomial  (Theorems \ref{thm:312}, \ref{thm:313}, \ref{thm:314} and \ref{thm:315}) are built with the help of the argument procedure of Theorems \ref{thm:34} and \ref{thm:38}.

\section{Copositivity of 4th order 2-dimensional symmetric tensors}

Let $T=(t_{i_1\cdots i_m}) (i_j=1,2,\ldots,n, j=1,2,\ldots, m)$ be a $m$th-order $n$-dimensional symmetric tensor.  Then for $x=(x_1.x_2,\cdots,x_n)^\top\in\mathbb{R}^n,$ we write
$$Tx^m=x^\top(Tx^{m-1})=\sum_{i_1\cdots i_m=1}^nt_{i_1\cdots i_m}x_{i_1}\cdots x_{i_m},$$
and $Tx^{m-1}=(y_1,y_2,\cdots, y_n)^\top$ is a vector with its components
$$y_k=(Tx^{m-1})_k=\sum_{i_2\cdots i_m=1}^nt_{ki_2\cdots i_m} x_{i_2}\cdots x_{i_m}, \  k=1,2,\cdots,n.$$

Let $f(x_1,x_2)$ be a quartic homogeneous real polynomial about two variables $x_1,x_2$,
\begin{equation}\label{eq:f}
	f(x_1,x_2)= x_1^4+4a x_1^3x_2+6b x_1^2x_2^2+4cx_1x_2^3+x_2^4.
\end{equation}
Then it gives a 4th-order 2-dimensional symmetric tensor $\mathcal{T}=(t_{ijkl})$ with its entiries,
\begin{equation}\label{eq:T2}
	t_{1111}=1, 
	t_{1112}=a, t_{1122}=b, t_{1222}=c, t_{2222}=1.
\end{equation}
By Theorem \ref{thm:13}, the following lemma can be obtained easily.
\begin{lemma}\label{lem:21}  Let $\mathcal{T}=(t_{ijkl})$  be a 4th-order 2-dimensional symmetric tensor given by \eqref{eq:T2}. Then $\mathcal{T}$ is copositive, i.e., $f(x_1,x_2)\geq0$ for all $x=(x_1,x_2)^\top\geq0$ if and only if
	\begin{itemize}
		\item[(1)] $\Delta'\leq0$ and $a+c>0;$  
		\item[(2)] $a\geq0$, $c\geq0$ and $1+3b\geq0;$  
		\item[(3)] $\Delta'\geq0$, $|a-c|\leq\sqrt{6b+2}$ and
		(i) $-1\leq 3b \leq3$, (ii) $b>1$ and $a+c \geq -\sqrt{6b-2}$,
	\end{itemize}
	where $\Delta=4\times12^3\Delta',$   $\Delta'=(1-4ac+3b^{2})^{3}-27(b+2abc-b^{3}-c^{2}-a^{2})^{2}$.
\end{lemma}

\begin{lemma}\label{lem:22} Let $\mathcal{T}=(t_{ijkl})$  be a 4th-order 2-dimensional symmetric tensor with its entires $|t_{ijkl}|=1$ and $t_{1111}=t_{2222}=1.$  Then $\mathcal{T}$ is copositive if and only if
	either $$b=t_{1122}=1\mbox{ or } a= t_{1112}=t_{1222}=c=1.$$\end{lemma}
{\bf Proof.} It follows from Lemma \ref{lem:21} that $\mathcal{T}$ is copositive if and only if
	\begin{itemize}
		\item[(1)] $\Delta'\leq0$ and $a+c>0;$  
		\item[(2)] $a\geq0$, $c\geq0$ and $1+3b\geq0 ;$  
		\item[(3)] $\Delta'\geq0$, $|a-c|\leq\sqrt{6b+2}$ and $-1\leq 3b \leq3;$.
	\end{itemize}
	
	Since  $|t_{ijkl}|=1$, then the conditions (1)-(3) mean 
	\begin{flushleft}
		\quad(1) $\Delta'\leq0$ and $a=c=1$ $\Leftrightarrow$  $a=c=1$ and either $$b=1, \Delta'=(1-4+3)^{3}-27(1+2-1^3-1-1)^{2}=0,$$  or $$b=-1, \Delta'=(1-4+3)^{3}-27(-1-2+1-1-1)^{2}<0;$$
	\end{flushleft}
	\begin{flushleft}
		\quad(2) $a=c=1$  and $b=1;$  
	\end{flushleft}
	\begin{flushleft}
		\quad(3) $\Delta'\geq0$, $|a-c|\leq\sqrt{6b+2}$ and $b=1$ $\Leftrightarrow$   $b=1$ and either $ac=1$,  $$ \Delta'=(1-4+3)^{3}-27(1+2-1-1-1)^{2}=0, |a-c|=0<\sqrt{6b+2}=\sqrt8;$$ or  $ac=-1,$
		$$ \Delta'=(1+4+3)^{3}-27(1-2-1-1-1)^{2}>0, |a-c|=2<\sqrt{6b+2}=\sqrt8.$$
	\end{flushleft}
	So the conditions (1)-(3) are equivalent to $$b=1\mbox{ or } a= c=1.$$
	This  completes the proof.	 

\begin{corollary}\label{cor:1} Let $\mathcal{T}=(t_{ijkl})$  be a 4th-order 2-dimensional symmetric tensor with its entires  $t_{1111}\geq1$ and $ t_{2222}\geq1.$  Then \begin{itemize}
		\item[(1)] $\mathcal{T}$ is strictly copositive if 
		$$ t_{1112}\geq1, t_{1222}\geq1, t_{1122}\geq-1;$$
		\item[(2)]	 $\mathcal{T}$ is copositive if 
		$$t_{1112}\geq-1, t_{1222}\geq-1, t_{1122}\geq1.$$ 
	\end{itemize}
\end{corollary}
{\bf Proof.}Let $x=(x_1,x_2)^\top\geq0.$ Then
	$$\mathcal{T}x^4= t_{1111}x_1^4+4t_{1112}x_1^3x_2+6t_{1122} x_1^2x_2^2+4t_{1222}x_1x_2^3+t_{2222}x_2^4.$$ 
	
	(1)	It is obvious that $$\begin{aligned}
		\mathcal{T}x^4\geq& x_1^4+4x_1^3x_2-6 x_1^2x_2^2+4x_1x_2^3+x_2^4
		=(x_1^2+x_2^2)^2+4x_1x_2(x_1^2-2x_1x_2+x_2^2)\\
		=&(x_1^2+x_2^2)^2+4x_1x_2(x_1-x_2)^2.
	\end{aligned}$$ So $\mathcal{T}x^4>0$ for $x\geq0$ and $x\ne0$. Suppose not, then there exists $x=(x_1,x_2)^\top\ne0$ such that $\mathcal{T}x^4=0,$ and hence, $$0=\mathcal{T}x^4\geq(x_1^2+x_2^2)^2+4x_1x_2(x_1-x_2)^2\geq0.$$That is, $$(x_1^2+x_2^2)^2+4x_1x_2(x_1-x_2)^2=0,$$ which means $x_1^2+x_2^2=0$, i.e., $x_1=x_2=0,$ a contradiction. Therefore, $\mathcal{T}$ is strictly copositive. 
	
	(2) For any $x\geq0$, it follows from Lemma \ref{lem:22} that  $$
	\mathcal{T}x^4\geq x_1^4-4x_1^3x_2+6 x_1^2x_2^2-4x_1x_2^3+x_2^4=(x_1-x_2)^4\geq0.$$ 
	So, $\mathcal{T}$ is copositive.
	This  completes the proof.	 

For a 4th-order 2-dimensional symmetric tensor $\mathcal{A}=(a_{ijkl})$   with its entires  $a_{1111}>0$ and $ a_{2222}>0,$ let $t_{1111}=t_{2222}=1$ and $$t_{1112}=a_{1112}a_{1111}^{-\frac34} a_{2222}^{-\frac14}, t_{1122}=a_{1122}a_{1111}^{-\frac12} a_{2222}^{-\frac12}, t_{1222}=a_{1222}a_{1111}^{-\frac14} a_{2222}^{-\frac34}.$$
For $y=(y_1,y_2)^\top$ and $x=(x_1,x_2)^\top=(a_{1111}^{\frac14}y_1,a_{2222}^{\frac14}y_2)^\top,$ then
$$\aligned \mathcal{A}y^4=& a_{1111}y_1^4+4a_{1112}y_1^3y_2+6a_{1122} y_1^2y_2^2+4a_{1222}y_1y_2^3+a_{2222}y_2^4\\
=&x_1^4+4t_{1112}x_1^3x_2+6t_{1122} x_1^2x_2^2+4t_{1222}x_1x_2^3+x_2^4=\mathcal{T}x^4.\endaligned$$ 
Obviously, the copositivity of $\mathcal{A}=(a_{ijkl})$ coincides with one of $\mathcal{T}=(t_{ijkl})$, and hence, we can establish an analytically sufficient condition of the (strict) copositivity of a gerenal 4th-order 2-dimensional symmetric tensor $\mathcal{A}=(a_{ijkl})$, that can be very easily parsed and validated.
\begin{corollary}\label{cor:2} Let $\mathcal{A}=(a_{ijkl})$  be a 4th-order 2-dimensional symmetric tensor with its entires  $a_{1111}>0$ and $ a_{2222}>0.$  Then \begin{itemize}
		\item[(1)] $\mathcal{A}$ is strictly copositive if 
		$$ a_{1112}\geq a_{1111}^{\frac34} a_{2222}^{\frac14}, a_{1122}\geq-\sqrt{a_{1111}a_{2222}}, a_{1222}\geq a_{1111}^{\frac14} a_{2222}^{\frac34};$$
		\item[(2)]	 $\mathcal{A}$ is copositive if 
		$$a_{1112}\geq-a_{1111}^{\frac34}a_{2222}^{\frac14}, a_{1122}\geq\sqrt{a_{1111}a_{2222}}, a_{1222}\geq-a_{1111}^{\frac14} a_{2222}^{\frac34}.$$ 
	\end{itemize}
\end{corollary}

\section{Copositivity of 4th order 3-dimensional symmetric tensors}

\subsection{Analytical expressions of strict copositivity}

\begin{theorem}\label{thm:31} \em Let $\mathcal{T}=(t_{ijkl})$ be a 4th-order $3$-dimensional symmetric tensor with $|t_{ijkl}|=t_{iiii}=t_{iiij}=1$ for all $i,j,k,l\in\{1,2,3\}$.  If there is at most one $-1$ in $\{t_{1123}, t_{1223}, t_{1233}\}$, then  $\mathcal{T}$ is strictly copositive.
\end{theorem}
{\bf Proof.} Without loss the generality, let $t_{1123}=-1, t_{1223}=t_{1233}=1$. Then $$\aligned 
	\mathcal{T}x^4=&x_1^4+x_2^4+x_3^4+6t_{1122}x_1^2x_2^2+6t_{1133}x_1^2x_3^2+6t_{2233}x_2^2x_3^2+4x_1^3x_2+4x_1^3x_3+4x_1x_2^3+4x_1x_3^3\\
	&+4x_2^3x_3+4x_2x_3^3-12x_1^2x_2x_3+12x_1x_2^2x_3+12x_1x_2x_3^2,
	\endaligned$$  and so, $$
	\mathcal{T}x^3=\dfrac14\nabla\mathcal{T}x^4= \left(\aligned \sum_{j,k,l=1}^3t_{1jkl}x_jx_kx_l\\
	\sum_{j,k,l=1}^3t_{2jkl}x_jx_kx_l\\ \sum_{j,k,l=1}^3t_{3jkl}x_jx_kx_l
	\endaligned\right)$$
	It follows from Theorem \ref{thm:14} that we need only show that $\mathcal{T}=(t_{ijkl})$ is strictly semi-positive, i.e., for $x=(x_1,x_2,x_3)^\top\geq0$, there exists $k\in\{1,2,3\}$ such that $$x_k>0\mbox{ and }(\mathcal{T}x^3)_k>0.$$
	For each $i\in\{1,2,3\}$, we have
	$$\aligned (\mathcal{T}x^3)_1 =&\sum_{j,k,l=1}^3t_{1jkl}x_jx_kx_l=x_1^3+x_2^3+x_3^3+3t_{1122}x_1x_2^2+3t_{1133}x_1x_3^2+3x_1^2x_2+3x_1^2x_3\\&+3x_2^2x_3+3x_2x_3^2-6x_1x_2x_3\\
	\geq&x_1^3+x_2^3+x_3^3-3x_1x_2^2-3x_1x_3^2+3x_1^2x_2+3x_1^2x_3+3x_2^2x_3+3x_2x_3^2-6x_1x_2x_3\\
	=&(x_2+x_3-x_1)^3+2x_1^3=x_1^3+(x_2+x_3-x_1+x_1)((x_2+x_3-x_1)^2-x_1(x_2+x_3-x_1)+x_1^2)\\
	=&x_1^3+(x_2+x_3)((x_2+x_3-x_1)^2-x_1(x_2+x_3-x_1)+x_1^2);
	\endaligned$$
	
	$$\aligned(\mathcal{T}x^3)_2=&\sum_{j,k,l=1}^3t_{2jkl}x_jx_kx_l= x_1^3+x_2^3+x_3^3+3t_{1122}x_1^2x_2+3t_{2233}x_2x_3^2+3x_1x_2^2+3x_2^2x_3\\&-3x_1^2x_3+3x_1x_3^2+6x_1x_2x_3\\
	\geq&x_1^3+x_2^3+x_3^3-3x_1^2x_2-3x_2x_3^2+3x_1x_2^2+3x_2^2x_3-3x_1^2x_3+3x_1x_3^2+6x_1x_2x_3\\
	=&(x_1+x_3-x_2)^3+2x_2^3+12x_1x_2x_3-6x_1^2x_3
	=((x_1+x_3-x_2)^3+2x_2^3)+6x_1x_3(2x_2-x_1);
	\endaligned$$
	$$\aligned (\mathcal{T}x^3)_3=&\sum_{j,k,l=1}^3t_{3jkl}x_jx_kx_l= x_1^3+x_2^3+x_3^3+3t_{1133}x_1^2x_3+3t_{2233}x_2^2x_3+3x_1x_3^2+3x_2x_3^2\\&-3x_1^2x_2+3x_1x_2^2+6x_1x_2x_3\\
	\geq&x_1^3+x_2^3+x_3^3-3x_1^2x_3-3x_2^2x_3+3x_1x_3^2+3x_2x_3^2-3x_1^2x_2+3x_1x_2^2+6x_1x_2x_3\\
	=&(x_1+x_2-x_3)^3+2x_3^3+12x_1x_2x_3-6x_1^2x_2
	=((x_1+x_2-x_3)^3+2x_3^3)+6x_1x_2(2x_3-x_1).
	\endaligned$$
	So,  it follows that \begin{itemize}
		\item $x_1>0$, $(\mathcal{T}x^3)_1>0,$ which is done; otherwise,
		\item $x_1=0, x_2>0$, $(\mathcal{T}x^3)_2>0;$
		\item $x_1=0, x_3>0$, $(\mathcal{T}x^3)_3>0,$
	\end{itemize}
	and hence, $\mathcal{T}$ is strictly copositive.
 
We definite $$
\mathcal{T'}=(t'_{ijkl})\leq\mathcal{T}=(t_{ijkl}) \Leftrightarrow t'_{ijkl}\leq t_{ijkl}, \mbox{ for all }i,j,k,l.
$$ Then for all $x\in\mathbb{R}^n_+,$ we have
$$\mathcal{T'}x^4=\sum_{i,j,k,l=1}^nt'_{ijkl}x_ix_jx_kx_l\leq\sum_{i,j,k,l=1}^nt_{ijkl}x_ix_jx_kx_l=\mathcal{T}x^4.$$ So, the (strict) copositivity of a tensor $\mathcal{T'}$ implies  one of  $\mathcal{T}$. Therefore, from Theorem \ref{thm:31},  the following conclusions are obvious.

\begin{corollary}\label{cor:32} \em Let $\mathcal{T}=(t_{ijkl})$ be a 4th-order $3$-dimensional symmetric tensor.  If $t_{iiii}\geq1, t_{iiij}\geq1, t_{iijj}\geq-1 \mbox { for all } i,j\in\{1,2,3\}, i\ne j$  and  one of the following conditions,  \begin{itemize}
		\item [(1)] $t_{1123}\geq-1, t_{1223}\geq1, t_{1233}\geq1$;
		\item [(2)]$t_{1123}\geq1, t_{1223}\geq-1, t_{1233}\geq1$; 
		\item [(3)] $t_{1123}\geq1, t_{1223}\geq1, t_{1233}\geq-1$,
	\end{itemize} 
	then  $\mathcal{T}$ is strictly copositive.
\end{corollary}

\begin{theorem}\label{thm:33} \em Let $\mathcal{T}=(t_{ijkl})$ be a 4th-order $3$-dimensional symmetric tensor with $|t_{ijkl}|=t_{iiii}=t_{iiij}=1$ for all $i,j,k,l\in\{1,2,3\}$. Then  $\mathcal{T}$ is srictly copositive  if and only if 
	\begin{itemize}
		\item [(1)] $t_{1122}=t_{1133}=t_{2233}=1 $  if $t_{1123}=t_{1223}=t_{1233}=-1$;
		\item [(2)] there is at least one $1$  in $\{t_{1122}, t_{1133}, t_{2233}\}$ and $t_{iijk}= t_{ijjk}=-t_{iijj}=-1$ for $i\ne j$ and $k\ne j$ and $i\ne k$   if two of  $\{t_{1123}, t_{1223}, t_{1233}\}$ are only $-1$.
	\end{itemize}
\end{theorem}

{\bf Proof.} Necessity. If $\mathcal{T}$ is srictly copositive, but the conditions don't hold. Then  

 (1)  for $x=(1,1,1)^\top$, it follows that    one of  $\{t_{1122}, t_{1133}, t_{2233}\}$ is $-1 $, $$\aligned 
	\mathcal{T}x^4=&x_1^4+x_2^4+x_3^4+6t_{1122}x_1^2x_2^2+6t_{1133}x_1^2x_3^2+6t_{2233}x_2^2x_3^2+4x_1^3x_2+4x_1^3x_3+4x_1x_2^3+4x_1x_3^3\\
	&+4x_2^3x_3+4x_2x_3^3+12t_{1123}x_1^2x_2x_3+12t_{1223}x_1x_2^2x_3+12t_{1233}x_1x_2x_3^2\\=&27+12-6-36=-3<0;
	\endaligned$$
	(2) for $x=(1,1,1)^\top$, it follows that  $t_{1122}= t_{1133}=t_{2233}=-1 $, $$\aligned 
	\mathcal{T}x^4=&x_1^4+x_2^4+x_3^4+6t_{1122}x_1^2x_2^2+6t_{1133}x_1^2x_3^2+6t_{2233}x_2^2x_3^2+4x_1^3x_2+4x_1^3x_3+4x_1x_2^3+4x_1x_3^3\\
	&+4x_2^3x_3+4x_2x_3^3+12t_{1123}x_1^2x_2x_3+12t_{1223}x_1x_2^2x_3+12t_{1233}x_1x_2x_3^2\\=&27-18-24+12=-3<0;
	\endaligned$$
	 Let $t_{1133}=1, t_{1122}=t_{2233}=-1$, $t_{1233}=1, t_{1223}=t_{1123}=-1$.  For  $x=(2,2.1,1)^\top$,
	$$\aligned 
	\mathcal{T}x^4=&(x_1+x_2+x_3)^4-12(x_1^2x_2^2+x_2^2x_3^2)-24(x_1^2x_2x_3+x_1x_2^2x_3)\\
	=&(2+2.1+1)^4-12(2^2\times 2.1^2+2.1^2)-24(2^2\times 2.1+2\times 2.1^2)\\
	=&-1.3599<0.
	\endaligned$$
	So, $\mathcal{T}$ is not srictly copositive.
	
	Sufficiency. (1) Since $t_{1122}=t_{1133}=t_{2233}=1 $ and  $t_{1123}=t_{1223}=t_{1233}=-1$, then $$\aligned 
	\mathcal{T}x^4=&x_1^4+x_2^4+x_3^4+6x_1^2x_2^2+6x_1^2x_3^2+6x_2^2x_3^2+4x_1^3x_2+4x_1^3x_3+4x_1x_2^3+4x_1x_3^3+4x_2^3x_3+4x_2x_3^3\\&-12x_1^2x_2x_3-12x_1x_2^2x_3-12x_1x_2x_3^2,
	\endaligned$$  and so, 
	$$\aligned (\mathcal{T}x^3)_1=&\sum_{j,k,l=1}^3t_{1jkl}x_jx_kx_l= x_1^3+x_2^3+x_3^3+3x_1x_2^2+3x_1x_3^2+3x_1^2x_2+3x_1^2x_3\\&-3x_2^2x_3-3x_2x_3^2-6x_1x_2x_3\\
	=&(x_2+x_3-x_1)^3+2x_1^3+6x_1x_2^2+6x_1x_3^2-6x_2^2x_3-6x_2x_3^2\\=&(x_2+x_3-x_1)^3+2x_1^3+6x_2^2(x_1-x_3)+6x_3^2(x_1-x_2);
	\endaligned$$
	$$\aligned(\mathcal{T}x^3)_2=&\sum_{j,k,l=1}^3t_{2jkl}x_jx_kx_l
	= x_1^3+x_2^3+x_3^3+3x_1^2x_2+3x_2x_3^2+3x_1x_2^2+3x_2^2x_3\\&-3x_1^2x_3-6x_1x_2x_3-3x_1x_3^2\\
	=&(x_1+x_3-x_2)^3+2x_2^3+6x_1^2x_2+6x_2x_3^2-6x_1^2x_3-6x_1x_3^2\\
	=&(x_1+x_3-x_2)^3+2x_2^3+6x_1^2(x_2-x_3)+6x_3^2(x_2-x_1);
	\endaligned$$
	$$\aligned(\mathcal{T}x^3)_3=&\sum_{j,k,l=1}^3t_{3jkl}x_jx_kx_l=x_1^3+x_2^3+x_3^3+3x_1^2x_3+3x_2^2x_3+3x_1x_3^2\\&+3x_2x_3^2-3x_1^2x_2-3x_1x_2^2-6x_1x_2x_3\\
	=&(x_1+x_2-x_3)^3+2x_3^3+6x_1^2x_3+6x_2^2x_3-6x_1^2x_2-6x_1x_2^2\\
	=&(x_1+x_2-x_3)^3+2x_3^3+6x_1^2(x_3-x_2)+6x_2^2(x_3-x_1).
	\endaligned$$
	So,  it follows that \begin{itemize}
		\item $x_1\geq \max\{x_2,x_3\}$, $(\mathcal{T}x^3)_1>0,$ which is done; otherwise,
		\item  $x_2\geq \max\{x_1,x_3\}$, $(\mathcal{T}x^3)_2>0;$
		\item  $x_3\geq \max\{x_1,x_2\}$, $(\mathcal{T}x^3)_3>0;$
	\end{itemize}
	and hence, $\mathcal{T}$ is strictly copositive. 
	
	(2) We might take $t_{1122}=t_{1133}=t_{1233}=1 $ and $t_{2233}=t_{1123}=t_{1223}=-1$, then $$\aligned 
	\mathcal{T}x^4=&x_1^4+x_2^4+x_3^4+6x_1^2x_2^2+6x_1^2x_3^2-6x_2^2x_3^2+4x_1^3x_2+4x_1^3x_3+4x_1x_2^3+4x_1x_3^3+4x_2^3x_3+4x_2x_3^3\\&-12x_1^2x_2x_3-12x_1x_2^2x_3+12x_1x_2x_3^2\\
	\geq &x_1^4+x_2^4+x_3^4+6x_1^2x_2^2-6x_1^2x_3^2-6x_2^2x_3^2+4x_1^3x_2+4x_1^3x_3+4x_1x_2^3+4x_1x_3^3+4x_2^3x_3+4x_2x_3^3\\&-12x_1^2x_2x_3-12x_1x_2^2x_3+12x_1x_2x_3^2=	\mathcal{T}'x^4,
	\endaligned$$  and so, 
	$$\aligned (\mathcal{T}'x^3)_1=&\sum_{j,k,l=1}^3t_{1jkl}x_jx_kx_l= x_1^3+x_2^3+x_3^3+3x_1x_2^2-3x_1x_3^2+3x_1^2x_2+3x_1^2x_3\\&-3x_2^2x_3+3x_2x_3^2-6x_1x_2x_3\\
	=&(x_2+x_3-x_1)^3+2x_1^3+6x_1x_2^2-6x_2^2x_3\\=&(x_2+x_3-x_1)^3+2x_1^3+6x_2^2(x_1-x_3);
	\endaligned$$
	$$\aligned(\mathcal{T}'x^3)_2=&\sum_{j,k,l=1}^3t_{2jkl}x_jx_kx_l=x_1^3+x_2^3+x_3^3+3x_1^2x_2-3x_2x_3^2+3x_1x_2^2+3x_2^2x_3\\&-3x_1^2x_3-6x_1x_2x_3+3x_1x_3^2\\
	=&(x_1+x_3-x_2)^3+2x_2^3+6x_1^2x_2-6x_1^2x_3\\
	=&(x_1+x_3-x_2)^3+2x_2^3+6x_1^2(x_2-x_3);
	\endaligned$$
	$$\aligned(\mathcal{T}'x^3)_3=&\sum_{j,k,l=1}^3t_{3jkl}x_jx_kx_l= x_1^3+x_2^3+x_3^3-3x_1^2x_3-3x_2^2x_3+3x_1x_3^2+3x_2x_3^2\\&-3x_1^2x_2-3x_1x_2^2+6x_1x_2x_3\\
	=&(x_1+x_2-x_3)^3+2x_3^3+12x_1x_2x_3-6x_1^2x_2-6x_1x_2^2\\
	=&(x_1+x_2-x_3)^3+2x_3^3+6x_1x_2((x_3-x_2)+(x_3-x_1)).
	\endaligned$$
	So,  it follows that \begin{itemize}
		\item $x_1\geq x_3$ and $x_1>0$, $(\mathcal{T}'x^3)_1>0,$ which is done; otherwise,
		\item  $x_2\geq x_3$ and $x_2>0$, $(\mathcal{T}'x^3)_2>0;$
		\item  $x_3\geq \max\{x_1,x_2\}$, $(\mathcal{T}'x^3)_3>0;$
	\end{itemize}
	and hence, $\mathcal{T}$ is strictly copositive. 
 
Combing the conclusions of Theorems \ref{thm:31} and \ref{thm:33}, the main result is bulit in this subsection.
\begin{theorem}\label{thm:34} \em Let $\mathcal{T}=(t_{ijkl})$ be a 4th-order $3$-dimensional symmetric tensor with $|t_{ijkl}|=t_{iiii}=t_{iiij}=1$ for all $i,j,k,l\in\{1,2,3\}$. Then  $\mathcal{T}$ is strictly copositive  if and only if 
	\begin{itemize}\item [(1)] there is at most one $-1$ in $\{t_{1123}, t_{1223}, t_{1233}\}$;
		\item [(2)]  two of  $\{t_{1123}, t_{1223}, t_{1233}\}$ are only $-1$ and  $t_{iijk}= t_{ijjk}=-t_{iijj}=-1$ for $i\ne j$ and $k\ne j$ and $i\ne k$;
		\item [(3)] $t_{1123}=t_{1223}=t_{1233}=-1$ and $t_{1122}=t_{1133}=t_{2233}=1 $.
	\end{itemize}
\end{theorem}

\subsection{Analytical expressions of copositivity}

\begin{theorem}\label{thm:35} \em Let $\mathcal{T}=(t_{ijkl})$ be a 4th-order $3$-dimensional symmetric tensor with $|t_{ijkl}|=t_{iiii}=t_{iijj}=1$ for all $i,j,k,l\in\{1,2,3\}$.  If $t_{1123}=t_{1223}=t_{1233}=1$, then  $\mathcal{T}$ is  copositive.
\end{theorem}
{\bf Proof.} Without loss the generality, let $x=(x_1,x_2,x_3)^\top\in\mathbb{R}^3_+$ with $x_1+x_2+x_3=1$. Then $$\aligned 
	\mathcal{T}x^4=&x_1^4+x_2^4+x_3^4+6x_1^2x_2^2+6x_1^2x_3^2+6x_2^2x_3^2+12x_1^2x_2x_3+12x_1x_2^2x_3+12x_1x_2x_3^2\\
	&+4t_{1112}x_1^3x_2+4t_{1113}x_1^3x_3+4t_{1222}x_1x_2^3+4t_{1333}x_1x_3^3+4t_{2223}x_2^3x_3+4t_{2333}x_2x_3^3,
	\endaligned$$  and so, 
	$$\aligned 
	\mathcal{T}x^4\geq&x_1^4+x_2^4+x_3^4+6x_1^2x_2^2+6x_1^2x_3^2+6x_2^2x_3^2-4x_1^3x_2-4x_1^3x_3-4x_1x_2^3-4x_1x_3^3\\&-4x_2^3x_3-4x_2x_3^3+12x_1^2x_2x_3+12x_1x_2^2x_3+12x_1x_2x_3^2\\
	=&(x_1+x_2-x_3)^4+8(3x_1^2x_2x_3+3x_1x_2^2x_3-x_1^3x_2-x_1x_2^3)\\
	=&(x_1+x_2+x_3)^4-8(x_1^3x_2+x_1x_2^3+x_1^3x_3+x_1x_3^3+x_2^3x_3+x_2x_3^3).
	\endaligned$$
	Let $f(x_1,x_2,x_3 )=(x_1+x_2-x_3)^4+8(3x_1^2x_2x_3+3x_1x_2^2x_3-x_1^3x_2-x_1x_2^3).$ 
	Solve the constrained  optimization problem in the non-negative orthant $\mathbb{R}^3_+$:
	$$\aligned
	\min \ &\ f(x_1,x_2,x_3)\\
	\mbox{ s. t. } & x_1+x_2+x_3=1.
	\endaligned$$
	Then the function $f(x_1,x_2,x_3)$ reaches  the minimum value $0$ at a point $\left(\dfrac12,0,\dfrac12\right)$ or $\left(0,\dfrac12,\dfrac12\right)$ or $\left(\dfrac12,\dfrac12,0\right)$, and  hence, $\mathcal{T}x^4\geq f(x_1,x_2,x_3)\geq0$ for any $x\geq0.$ That is, $\mathcal{T}$ is copositive.

From the proof of Theorem \ref{thm:35}, the following conclusion is easily obtained.
\begin{corollary}\label{cor:36} \em Let $\mathcal{T}=(t_{ijkl})$ be a 4th-order $3$-dimensional symmetric tensor.  If for all $i,j\in\{1,2,3\}$ and $i\ne j$, $$t_{iiii}\geq1, t_{iijj}\geq1, t_{iiij}\geq-1, t_{1123}\geq1, t_{1223}\geq1, t_{1233}\geq1$$
	then  $\mathcal{T}$ is copositive.
\end{corollary}

\begin{theorem}\label{thm:37} \em Let $\mathcal{T}=(t_{ijkl})$ be a 4th-order $3$-dimensional symmetric tensor with $|t_{ijkl}|=t_{iiii}=t_{iijj}=1$ for all $i,j,k,l\in\{1,2,3\}$. Then  $\mathcal{T}$ is copositive  if and only if \begin{itemize}
		\item [(1)] there is at least  one  $1$  in $\{t_{iiij}; i,j=1,2,3,i\ne j\}$ and $t_{iijk}=t_{iiij}t_{iiik}=-1$ or $t_{iijk}=-t_{iiij}=-t_{iiik}=-1$for $i\ne j$ and $k\ne j$ and $i\ne k$ if there is only one $-1$  in $\{t_{1123}, t_{1223}, t_{1233}\}$;
		\item [(2)] there is at least  two  $1$  in $\{t_{iiij}; i,j=1,2,3,i\ne j\}$ and there is at least  one  $1$  in $\{t_{iiij},t_{iiik};t_{iijk}=-1,i\ne j,i\ne k, i\ne k\}$  if there is only two  $-1$  in $\{t_{1123}, t_{1223}, t_{1233}\}$;
		\item [(3)]  $t_{iiij}=1$ for all $i,j=1,2,3$ with $i\ne j$  if $t_{1123}=t_{1223}=t_{1233}=-1$.
	\end{itemize}
\end{theorem}
{\bf Proof.}  Necessity. Suppose $\mathcal{T}$ is copositive, but the conditions don't hold. 
	
	(1)  Assume $t_{iiij}=-1$ for all $ i,j=1,2,3, i\ne j$. Then  we might take $t_{1123}=-1$,  and so for $x=(3,1,1)^\top$, $$\aligned 
	\mathcal{T}x^4=&x_1^4+x_2^4+x_3^4+6x_1^2x_2^2+6x_1^2x_3^2+6x_2^2x_3^2-4x_1^3x_2-4x_1^3x_3-4x_1x_2^3-4x_1x_3^3\\&-4x_2^3x_3-4x_2x_3^3-12x_1^2x_2x_3+12x_1x_2^2x_3+12x_1x_2x_3^2\\
=&(x_1+x_2-x_3)^4+8(3x_1x_2^2x_3-x_1^3x_2-x_1x_2^3)\\
	=&(3+1-1)^4+8(3^2-3^3-3)<0;
	\endaligned$$
	Assume $t_{iijk}=t_{iiij}=t_{iiik}=-1$.	 We might take $t_{1113}=t_{1123}=t_{1233}=1$, $t_{1223}=t_{1222}=t_{2223}=t_{1333}=t_{1112}=t_{2333}=-1$, then for $x=(1,3,1)^\top$,
	$$\aligned 
	\mathcal{T}x^4=&x_1^4+x_2^4+x_3^4+6x_1^2x_2^2+6x_1^2x_3^2+6x_2^2x_3^2-4x_1^3x_2+4x_1^3x_3-4x_1x_2^3-4x_1x_3^3\\&-4x_2^3x_3-4x_2x_3^3+12x_1^2x_2x_3-12x_1x_2^2x_3+12x_1x_2x_3^2\\
	=&(x_1+x_2-x_3)^4+8(3x_1^2x_2x_3-x_1^3x_2-x_1x_2^3+x_1^3x_3)\\
	=&(1+3-1)^4+8(3^2-3-3^3+1)<0.\endaligned$$
	(2) If  there are five  $-1$  in $\{t_{iiij}; i,j=1,2,3,i\ne j\}$,  then for $x=(1,1,1)^\top$,  $$\aligned 
	\mathcal{T}x^4=&x_1^4+x_2^4+x_3^4+6x_1^2x_2^2+6x_1^2x_3^2+6x_2^2x_3^2\\&+4t_{1112}x_1^3x_2+4t_{1113}x_1^3x_3+4t_{1222}x_1x_2^3+4t_{1333}x_1x_3^3+4t_{2223}x_2^3x_3+4t_{2333}x_2x_3^3\\&+12t_{1123}x_1^2x_2x_3+12t_{1223}x_1x_2^2x_3+12t_{1233}x_1x_2x_3^2\\=&21+4-20-24+12=-7<0;
	\endaligned$$
	
Assume $t_{iiij}=t_{iiik}=t_{iijk}=-1.$ We might take $t_{1333}=t_{2333}=t_{1123}= 1$ and $t_{1113}=t_{1112}=t_{2223}=t_{1222}=t_{1233}=t_{1223}= -1$.  Then for $x=(1,3,1)^\top$,  $$\aligned 
	\mathcal{T}x^4=&(x_1+x_2+x_3)^4-8(3x_1x_2x_3^2+3x_1x_2^2x_3+x_1^3x_2+x_1x_2^3+x_1x_3^3+x_2^3x_3)\\
	=&5^4-8(9+27+3+27+1+27)=-127<0;
	\endaligned$$
	
	(3) Assume there is one  $-1$  in $\{t_{iiij}; i,j=1,2,3,i\ne j\}$, we might let $t_{1112}=-1$. Then for $x=(4,3,2)^\top$, 
$$\aligned \mathcal{T}x^4=&x_1^4+x_2^4+x_3^4+6x_1^2x_2^2+6x_1^2x_3^2+6x_2^2x_3^2-4x_1^3x_2+4x_1^3x_3+4x_1x_2^3+4x_1x_3^3\\
&+4x_2^3x_3+4x_2x_3^3-12x_1^2x_2x_3-12x_1x_2^2x_3-12x_1x_2x_3^2\\
	=&(x_1+x_2+x_3)^4-8(x_1^3x_2+3x_1^2x_2x_3+3x_1x_2^2x_3+3x_1x_2x_3^2)\\
	=&(4+3+2)^4-8(4^3\times3 + 3 \times4^2\times3\times2 + 3 \times4\times3^2\times2 + 3\times4\times3\times2^2)<0.
	\endaligned$$
	So, $\mathcal{T}$ is not copositive, which is a contradiction, and hence,  the conditions hold.
	
	Sufficiency. 
	$\mathcal{T}x^4 $ may be rewritten as follows, $$\aligned 
	\mathcal{T}x^4=&(x_1+x_2-x_3)^4+4(t_{1112}-1)x_1^3x_2+4(t_{1222}-1)x_1x_2^3+4(t_{1333}+1)x_1x_3^3\\&+4(t_{1113}+1)x_1^3x_3+4(t_{2223}+1)x_2^3x_3+4(t_{2333}+1)x_2x_3^3\\
	&+12(t_{1123}+1)x_1^2x_2x_3+12(t_{1223}+1)x_1x_2^2x_3+12(t_{1233}-1)x_1x_2x_3^2;
	\endaligned$$$$\aligned 
	\mathcal{T}x^4=&(x_1-x_2+x_3)^4+4(t_{1112}+1)x_1^3x_2+4(t_{1222}+1)x_1x_2^3+4(t_{1333}-1)x_1x_3^3\\
	&+4(t_{1113}-1)x_1^3x_3+4(t_{2223}+1)x_2^3x_3+4(t_{2333}+1)x_2x_3^3\\
	&+12(t_{1123}+1)x_1^2x_2x_3+12(t_{1223}-1)x_1x_2^2x_3+12(t_{1233}+1)x_1x_2x_3^2;
	\endaligned$$$$\aligned 
	\mathcal{T}x^4=&(x_1-x_2-x_3)^4+4(t_{1112}+1)x_1^3x_2+4(t_{1222}+1)x_1x_2^3+4(t_{1333}+1)x_1x_3^3\\
	&+4(t_{1113}+1)x_1^3x_3+4(t_{2223}-1)x_2^3x_3+4(t_{2333}-1)x_2x_3^3\\
	&+12(t_{1123}-1)x_1^2x_2x_3+12(t_{1223}+1)x_1x_2^2x_3+12(t_{1233}+1)x_1x_2x_3^2.
	\endaligned$$
	
	Clearly, for the boundary points of the non-negative orthant,  $x=(0,x_2,x_3)$, $(x_1,0,x_3)$, $(x_1,x_2,0)$,  it follows from Lemma \ref{lem:22} that $\mathcal{T}x^4\geq0$. Let $x_1+x_2+x_3=1$ in the sequal.
	
	(1) Without loss the generality, let $t_{1112}=t_{1233}=t_{1223}= 1$ and $t_{1113}=t_{1333}=t_{2333}=t_{1222}=t_{2223}=t_{1123}= -1.$  Then  
	$$ 
	\mathcal{T}x^4=(x_1+x_2-x_3)^4+8(3x_1x_2^2x_3-x_1x_2^3).
	$$
	Solve the constrained  optimization problem in the non-negative orthant $\mathbb{R}^3_+$:
	$$\aligned
	\min \ &\ \mathcal{T}x^4\\
	\mbox{ s. t. } & x_1+x_2+x_3=1.
	\endaligned$$
	Then the polynomial $\mathcal{T}x^4$ reaches  the minimum value $0$ at a point $\left(\dfrac12,0,\dfrac12\right)$ or $\left(0, \dfrac12,\dfrac12\right)$ , and  hence, $\mathcal{T}x^4\geq0$ for all  $x\geq0.$ That is, $\mathcal{T}$ is copositive.
	
	(2) Without loss the generality, let $t_{1112}=t_{2223}=t_{1233}= 1$ and $t_{1113}=t_{1333}=t_{2333}=t_{1222}=t_{1123}=t_{1223}= -1.$  Then 
	$$ 
	\mathcal{T}x^4=(x_1+x_2-x_3)^4+8x_2^3(x_3-x_1).
	$$
	Let $f(x_1,x_2,x_3,\lambda)=\mathcal{T}x^4+\lambda(x_1+x_2+x_3-1).$ Then the stationary points of the function $f(x_1,x_2,x_3,\lambda)$ are the solution  to this system of equations,
	$$\begin{cases}
		f'_{x_1}(x_1,x_2,x_3)=4(x_1+x_2-x_3)^3-8x_2^3+\lambda=0,\\
		f'_{x_2}(x_1,x_2,x_3)=4(x_1+x_2-x_3)^3+24x_2^2(x_3-x_1)+\lambda=0,\\
		f'_{x_3}(x_1,x_2,x_3)=-4(x_1+x_2-x_3)^3+8x_2^3+\lambda=0,\\
		x_1+x_2+x_3=1
	\end{cases}$$
Solve it in the non-negative orthant $\mathbb{R}^3_+$,  $$x_1=x_3=\dfrac12, x_2=0,\lambda=0.$$
	So $\mathcal{T}x^4\geq0$ at  the boundary points, and then  $\mathcal{T}x^4\geq0$ for all $x\geq0.$ That is, $\mathcal{T}$ is copositive.
	
	(3) It follows from There \ref{thm:33} (1) that $\mathcal{T}$ is copositive. This completes the proof.\\
 
Combing the conclusions of Theorems \ref{thm:35} and \ref{thm:37}, the main result is established in this subsection.
\begin{theorem}\label{thm:38} \em Let $\mathcal{T}=(t_{ijkl})$ be a 4th-order $3$-dimensional symmetric tensor with $|t_{ijkl}|=t_{iiii}=t_{iijj}=1$ for all $i,j,k,l\in\{1,2,3\}$. Then  $\mathcal{A}$ is copositive  if and only if \begin{itemize}
		\item[(1)] $t_{1123}=t_{1223}=t_{1233}=1$;
		\item [(2)] there is only one $-1$  in $\{t_{1123}, t_{1223}, t_{1233}\}$ and  $t_{iijk}=t_{iiij}t_{iiik}=-1$ or $t_{iijk}=-t_{iiij}=-t_{iiik}=-1$for $i\ne j$ and $k\ne j$ and $i\ne k$;
		\item [(3)]  there is only two  $-1$  in $\{t_{1123}, t_{1223}, t_{1233}\}$ and there is at least  one  $1$  in $\{t_{iiij},t_{iiik};t_{iijk}=-1,i\ne j,i\ne k, i\ne k\}$;
		\item [(4)] $t_{1123}=t_{1223}=t_{1233}=-1$ and  $t_{iiij}=1, i,j=1,2,3, i\ne j$.
	\end{itemize}
\end{theorem}

\subsection{Applications to a gerenal   tensor }
In this subsection, we apply Theorems \ref{thm:34} and \ref{thm:38} to find the (strict) copositivity of a gerenal 4th-order 3-dimensional symmetric tensor, and moreover, these analytic conditions can be very easily parsed and verified.

For a 4th-order 3-dimensional symmetric tensor $\mathcal{T}=(t_{ijkl})$   with its entires  $t_{iiii}>0$ for all $i\in\{1,2,3\},$ let  $\mathcal{T'}=(t'_{ijkl})$ be a symmetric tensor with its entires $t'_{1111}=t'_{2222}=t'_{3333}=1$ and $$t'_{1112}=t_{1112}t_{1111}^{-\frac34} t_{2222}^{-\frac14}, t'_{1122}=t_{1122}t_{1111}^{-\frac12} t_{2222}^{-\frac12}, t'_{1222}=t_{1222}t_{1111}^{-\frac14} t_{2222}^{-\frac34},$$
$$t'_{1113}=t_{1113}t_{1111}^{-\frac34} t_{3333}^{-\frac14}, t'_{1133}=t_{1133}t_{1111}^{-\frac12} t_{3333}^{-\frac12}, t'_{1333}=t_{1333}t_{1111}^{-\frac14} t_{3333}^{-\frac34},$$ 
$$t'_{2223}=t_{2223}t_{2222}^{-\frac34} t_{3333}^{-\frac14}, t'_{2233}=t_{2233}t_{2222}^{-\frac12} t_{3333}^{-\frac12}, t'_{2333}=t_{2333}t_{2222}^{-\frac14} t_{3333}^{-\frac34},$$
$$t'_{1123}=t_{1123}t_{1111}^{-\frac12} t_{2222}^{-\frac14} t_{3333}^{-\frac14}, t'_{1223}=t_{1223}t_{1111}^{-\frac14} t_{2222}^{-\frac12} t_{3333}^{-\frac14}, t'_{1233}=t_{1233}t_{1111}^{-\frac14} t_{2222}^{-\frac14} t_{3333}^{-\frac12}.$$
For $y=(y_1,y_2,y_3)^\top$ and $x=(x_1,x_2,x_3)^\top=(t_{1111}^{\frac14}y_1,t_{2222}^{\frac14}y_2,t_{3333}^{\frac14}y_3)^\top,$ then
$$\aligned \mathcal{T}y^4=& t_{1111}y_1^4+4t_{1112}y_1^3y_2+6t_{1122} y_1^2y_2^2+4t_{1222}y_1y_2^3+t_{2222}y_2^4+4t_{1113}y_1^3y_3+6t_{1133}y_1^2y_3^2\\&+4t_{1333} y_1y_3^3+t_{3333}y_3^4 +4t_{2223}y_2^3y_3+6t_{2233}y_2^2y_3^2+4t_{2333} y_2y_3^3 \\
& +12t_{1123}y_1^2y_2y_3+12t_{1223}y_1y_2^2y_3+12t_{1233} y_1y_2y_3^3\\
=&x_1^4+4t'_{1112}x_1^3x_2+6t'_{1122} x_1^2x_2^2+4t'_{1222}x_1x_2^3+x_2^4+4t'_{1113}x_1^3x_3+6t'_{1133} x_1^2x_3^2\\&+4t'_{1333}x_1x_3^3+x_3^4+4t'_{2223}x_2^3x_3+6t'_{2233} x_2^2x_3^2+4t'_{2333}x_2x_3^3\\
&+12t'_{1123}x_1^2x_2x_3+12t'_{1223} x_1x_2^2x_3+12t'_{1233}x_1x_2x_3^2\\
=&\mathcal{T'}x^4.\endaligned$$ 
It is obvious that the copositivity of symmetric  tensor $\mathcal{T}=(t_{ijkl})$ is equivalent to the copositivity of   $\mathcal{T'}=(t'_{ijkl})$.  So,  applying Corollaries \ref{cor:32} and \ref{cor:36} (or Theorems \ref{thm:34} (1) and \ref{thm:38}  (2) to establish easily the following conclusions .

\begin{theorem}\label{thm:39} \em Let $\mathcal{T}=(t_{ijkl})$ be a 4th-order $3$-dimensional symmetric tensor with $t_{iiii}>0$ for all $i\in\{1,2,3\}$ .  Assume one of the following three conditions holds,  \begin{itemize}
		\item [(a)] $t_{1123}\geq-t_{1111}^{\frac12}t_{2222}^{\frac14}t_{3333}^{\frac14}, t_{1223}\geq t_{1111}^{\frac14}t_{2222}^{\frac12}t_{3333}^{\frac14}, t_{1233}\geq t_{1111}^{\frac14}t_{2222}^{\frac14}t_{3333}^{\frac12}$; 
		\item [(b)]$t_{1123}\geq t_{1111}^{\frac12}t_{2222}^{\frac14}t_{3333}^{\frac14}, t_{1223}\geq-t_{1111}^{\frac14}t_{2222}^{\frac12}t_{3333}^{\frac14}, t_{1233}\geq t_{1111}^{\frac14}t_{2222}^{\frac14}t_{3333}^{\frac12}$; 
		\item [(c)] $t_{1123}\geq t_{1111}^{\frac12}t_{2222}^{\frac14}t_{3333}^{\frac14}, t_{1223}\geq t_{1111}^{\frac14}t_{2222}^{\frac12}t_{3333}^{\frac14}, t_{1233}\geq-t_{1111}^{\frac14}t_{2222}^{\frac14}t_{3333}^{\frac12}$.
	\end{itemize} 
	Then (1) $\mathcal{T}$ is strictly copositive if for all $i,j\in\{1,2,3\}$ and $i\ne j$, $$ t_{iijj}\geq -\sqrt{t_{iiii}t_{jjjj}}, \ \ t_{iiij}\geq t_{iiii}^{\frac34}t_{jjjj}^{\frac14};$$ 
	(2) $\mathcal{T}$ is copositive if  for all $i,j\in\{1,2,3\}$ and $i\ne j$, $$ t_{iijj}\geq \sqrt{t_{iiii}t_{jjjj}},\ \ t_{iiij}\geq-t_{iiii}^{\frac34}t_{jjjj}^{\frac14},$$ and there is at least one $t_{rrrs}\in\{t_{iiij}, t_{iiik}; t_{iijk}\geq-t_{iiii}^{\frac12}t_{jjjj}^{\frac14}t_{kkkk}^{\frac14}, i\ne j, i\ne k, k\ne j\}$ such that $$t_{rrrs}\geq t_{rrrr}^{\frac34}t_{ssss}^{\frac14}.$$
	\end{theorem}

From Theorems \ref{thm:34} (3) or \ref{thm:38} (4), the following conclusions are established easily.

\begin{theorem}\label{thm:310} \em Let $\mathcal{T}=(t_{ijkl})$ be a 4th-order $3$-dimensional symmetric tensor with $t_{iiii}>0$ for all $i\in\{1,2,3\}$ .  Assume that   $$t_{1123}\geq-t_{1111}^{\frac12}t_{2222}^{\frac14}t_{3333}^{\frac14}, t_{1223}\geq -t_{1111}^{\frac14}t_{2222}^{\frac12}t_{3333}^{\frac14}, t_{1233}\geq -t_{1111}^{\frac14}t_{2222}^{\frac14}t_{3333}^{\frac12}.$$
	Then   $\mathcal{T}$ is strictly copositive if for all $i,j\in\{1,2,3\}$ and $i\ne j$, $$ t_{iijj}\geq \sqrt{t_{iiii}t_{jjjj}}, \ \ t_{iiij}\geq t_{iiii}^{\frac34}t_{jjjj}^{\frac14}.$$
\end{theorem}

From Theorems \ref{thm:34} (2) and \ref{thm:38} (3), the following conclusions are established easily.

\begin{theorem}\label{thm:311} \em Let $\mathcal{T}=(t_{ijkl})$ be a 4th-order $3$-dimensional symmetric tensor with $t_{iiii}>0$ for all $i\in\{1,2,3\}$ .  Assume that   one of the following three conditions holds,  \begin{itemize}
		\item [(a)] $t_{1123}\geq-t_{1111}^{\frac12}t_{2222}^{\frac14}t_{3333}^{\frac14}, t_{1223}\geq -t_{1111}^{\frac14}t_{2222}^{\frac12}t_{3333}^{\frac14}, t_{1233}\geq t_{1111}^{\frac14}t_{2222}^{\frac14}t_{3333}^{\frac12}$; 
		\item [(b)]$t_{1123}\geq t_{1111}^{\frac12}t_{2222}^{\frac14}t_{3333}^{\frac14}, t_{1223}\geq-t_{1111}^{\frac14}t_{2222}^{\frac12}t_{3333}^{\frac14}, t_{1233}\geq- t_{1111}^{\frac14}t_{2222}^{\frac14}t_{3333}^{\frac12}$; 
		\item [(c)] $t_{1123}\geq -t_{1111}^{\frac12}t_{2222}^{\frac14}t_{3333}^{\frac14}, t_{1223}\geq t_{1111}^{\frac14}t_{2222}^{\frac12}t_{3333}^{\frac14}, t_{1233}\geq-t_{1111}^{\frac14}t_{2222}^{\frac14}t_{3333}^{\frac12}$.
	\end{itemize} 
	Then (1) $\mathcal{T}$ is strictly copositive if for all $i,j\in\{1,2,3\}$ and $i\ne j$, $$  t_{iiij}\geq t_{iiii}^{\frac34}t_{jjjj}^{\frac14};$$
	and $$t_{iijj}\geq\sqrt{t_{iiii}t_{jjjj}} \mbox{ and } t_{iijk}\geq-t_{iiii}^{\frac12}t_{jjjj}^{\frac14}t_{kkkk}^{\frac14}\mbox{ and } t_{ijjk}\geq -t_{iiii}^{\frac14}t_{jjjj}^{\frac12}t_{kkkk}^{\frac14}.$$ 
	(2) $\mathcal{T}$ is copositive if for all $i,j\in\{1,2,3\}$ and $i\ne j$, $$ t_{iijj}\geq \sqrt{t_{iiii}t_{jjjj}}$$ and there is at least one $t_{sssr}\in\{t_{iiij}, t_{iiik}; t_{iijk}\geq-t_{iiii}^{\frac12}t_{jjjj}^{\frac14}t_{kkkk}^{\frac14}, i\ne j, i\ne k, k\ne j\}$ such that $$t_{sssr}\geq t_{ssss}^{\frac34}t_{rrrr}^{\frac14}.$$
\end{theorem}
\subsection{Ternary quartic inequalities }

In this subsection, we apply Theorems \ref{thm:34} and \ref{thm:38} to show several ternary quartic (strict) inequalities.
Applying Theorems \ref{thm:34} to establish easily the following strict inequalities .

\begin{theorem}\label{thm:312}	If $(x_1,x_2,x_3)\ne(0,0,0)$ and $x_i\geq0,\ i=1,2,3$, then 
	\begin{itemize}
		\item [(i)] $x_1^4+x_2^4+x_3^4+6x_1^2x_2^2+6x_1^2x_3^2+6x_2^2x_3^2+
		4x_1^3x_3+4x_1x_2^3+4x_1x_3^3+4x_2^3x_3+4x_2x_3^3+4x_1^3x_2>12x_1^2x_2x_3+12x_1x_2^2x_3+12x_1x_2x_3^2,
		$ \ \  or equivalently,
		$$(x_1+x_2+x_3)^4>24x_1x_2x_3(x_1+x_2+x_3);$$
		\item [(ii)] $x_1^4+x_2^4+x_3^4+
		4x_1^3x_3+4x_1x_2^3+4x_1x_3^3+4x_2^3x_3+4x_2x_3^3+4x_1^3x_2+12x_1x_2^2x_3+12x_1x_2x_3^2>6x_1^2x_2^2+6x_1^2x_3^2+6x_2^2x_3^2+12x_1^2x_2x_3,$ \ \  or equivalently,
		$$(x_1+x_2+x_3)^4>12(x_1^2x_2^2+x_1^2x_3^2+x_2^2x_3^2+2x_1^2x_2x_3);$$
		\item [(iii)]$x_1^4+x_2^4+x_3^4+
		4x_1^3x_3+4x_1x_2^3+4x_1x_3^3+4x_2^3x_3+4x_2x_3^3+4x_1^3x_2+12x_1x_2^2x_3+12x_1^2x_2x_3>6x_1^2x_2^2+6x_1^2x_3^2+6x_2^2x_3^2+12x_1x_2x_3^2,
		$\ \  or equivalently,$$(x_1+x_2+x_3)^4>12(x_1^2x_2^2+x_1^2x_3^2+x_2^2x_3^2+2x_1x_2x_3^2);$$
		\item [(iv)] $x_1^4+x_2^4+x_3^4+
		4x_1^3x_3+4x_1x_2^3+4x_1x_3^3+4x_2^3x_3+4x_2x_3^3+4x_1^3x_2+12x_1^2x_2x_3+12x_1x_2x_3^2>6x_1^2x_2^2+6x_1^2x_3^2+6x_2^2x_3^2+12x_1x_2^2x_3,
		$\ \  or equivalently,$$(x_1+x_2+x_3)^4>12(x_1^2x_2^2+x_1^2x_3^2+x_2^2x_3^2+2x_1x_2^2x_3);$$
		\item [(v)] $x_1^4+x_2^4+x_3^4+6x_1^2x_3^2+
		4x_1^3x_3+4x_1x_2^3+4x_1x_3^3+4x_2^3x_3+4x_2x_3^3+4x_1^3x_2+12x_1x_2^2x_3>6x_1^2x_2^2+6x_2^2x_3^2+12x_1^2x_2x_3+12x_1x_2x_3^2,
		$ \ \  or equivalently,$$(x_1-x_2+x_3)^4>12(x_1^2x_2^2+x_2^2x_3^2);$$
		\item [(vi)]$x_1^4+x_2^4+x_3^4+6x_2^2x_3^2+
		4x_1^3x_3+4x_1x_2^3+4x_1x_3^3+4x_2^3x_3+4x_2x_3^3+4x_1^3x_2+12x_1^2x_2x_3>6x_1^2x_3^2+6x_1^2x_2^2+12x_1x_2^2x_3+12x_1x_2x_3^2,
		$\ \  or equivalently,$$(x_2-x_1+x_3)^4>12(x_1^2x_3^2+x_2^2x_1^2);$$
		\item [(vii)] $x_1^4+x_2^4+x_3^4+6x_1^2x_2^2+
		4x_1^3x_3+4x_1x_2^3+4x_1x_3^3+4x_2^3x_3+4x_2x_3^3+4x_1^3x_2+12x_1x_2x_3^2>6x_1^2x_3^2+6x_2^2x_3^2+12x_1x_2^2x_3+12x_1^2x_2x_3,
		$\ \  or equivalently,$$(x_1+x_2-x_3)^4>12(x_1^2x_3^2+x_2^2x_3^2).$$
	\end{itemize} 
\end{theorem}

It follows from Theorems  \ref{thm:38} that the following inequalities are obtained easily.
\begin{theorem}\label{thm:313} \em Let $x_1,x_2,x_3$ be three nonnegative real numbers.  Then   \begin{itemize}
		\item [(a)] $x_1^4+x_2^4+x_3^4+6x_1^2x_2^2+6x_1^2x_3^2+6x_2^2x_3^2+12x_1^2x_2x_3+12x_1x_2x_3^2+12x_1x_2^2x_3\geq
		4x_1^3x_3+4x_1x_2^3+4x_2^3x_3+4x_2x_3^3+4x_1^3x_2+4x_1x_3^3,
		$ \ \  or equivalently,$$(x_1+x_2+x_3)^4\geq8(x_1^3x_3 +x_1x_2^3+x_2^3x_3+x_2x_3^3+x_1^3x_2+x_1x_3^3)$$
		with equality if and only if  two of $x_1, x_2, x_3$	are equal and the third is $0$.
		\item [(b)] $x_1^4+x_2^4+x_3^4+6x_1^2x_2^2+6x_1^2x_3^2+6x_2^2x_3^2+12x_1^2x_2x_3+12x_1x_2x_3^2+4x_1x_2^3\geq
		4x_1^3x_3+4x_1x_3^3+4x_2^3x_3+4x_2x_3^3+4x_1^3x_2+12x_1x_2^2x_3,
		$ \ \  or equivalently,$$(x_1+x_2-x_3)^4\geq8x_1^2x_2(x_1-3 x_3)$$
		with equality if and only if  $x_1= 0$, $x_3=x_2$ or $x_2=0,$ $x_1=x_3$.
				\item [(c)] $x_1^4+x_2^4+x_3^4+6x_1^2x_2^2+6x_1^2x_3^2+6x_2^2x_3^2+12x_1x_2^2x_3+12x_1x_2x_3^2+4x_1^3x_3\geq
		4x_1x_2^3+4x_1x_3^3+4x_2^3x_3+4x_2x_3^3+4x_1^3x_2+12x_1^2x_2x_3,
		$ \ \  or equivalently,$$(x_1-x_2+x_3)^4\geq8x_1x_3^2( x_3-3x_2)$$
		with equality if and only if $x_1= x_2$, $x_3=0$ or $x_1=0,$ $x_2=x_3$.
		\item [(d)]$x_1^4+x_2^4+x_3^4+6x_1^2x_2^2+6x_1^2x_3^2+6x_2^2x_3^2+12x_1x_2^2x_3+12x_1^2x_2x_3+4x_2x_3^3\geq
		4x_1^3x_3+4x_1x_3^3+4x_2^3x_3+4x_2x_3^3+4x_1^3x_2+12x_1x_2x_3^2,
		$ \ \  or equivalently,$$(x_2+x_3-x_1)^4\geq8x_2^2x_3(x_2-3x_1)$$
		with equality if and only if  $x_1= x_2$, $x_3=0$ or $x_2=0,$ $x_1=x_3$.
			\end{itemize} 
	
\end{theorem}

{\bf Proof.} It follows from Theorems  \ref{thm:38} (1) and  (2) that the  inequalities hold easily. Now we show the equality of (a), and the arguments of (b), (c) and (d) are the same. 	$$\aligned 
	f(x_1,x_2,x_3)=&x_1^4+x_2^4+x_3^4+6x_1^2x_2^2+6x_1^2x_3^2+6x_2^2x_3^2-4x_1^3x_2-4x_1^3x_3-4x_1x_2^3-4x_1x_3^3\\&-4x_2^3x_3-4x_2x_3^3+12x_1^2x_2x_3+12x_1x_2^2x_3+12x_1x_2x_3^2\\
	=&(x_1+x_2+x_3)^4-8(x_1^3x_3 +x_1x_2^3+x_2^3x_3+x_2x_3^3+x_1^3x_2+x_1x_3^3).
	\endaligned$$
	From the arguments of Theorems  \ref{thm:35}, it follows that the zero of function $f(x_1,x_2,x_3)$  must be on the boundary of the non-negative orthant (that is, three coordinate planes), and so, $f(x_1,x_2,x_3)=0$ if and only if  $x_1= x_2$, $x_3=0$ or $x_1=0,$ $x_2=x_3$ or $x_1=x_3, x_2=0$.\\

Similarly, from the arguments of Theorems  \ref{thm:37} (2), we  have the following conclusion also.
\begin{theorem}\label{thm:314} \em Let $x_1,x_2,x_3$ be three nonnegative real numbers.  Then   \begin{itemize}
		\item [(e)] $x_1^4+x_2^4+x_3^4+6x_1^2x_2^2+6x_1^2x_3^2+6x_2^2x_3^2+4x_1^3x_2+4x_2^3x_3+12x_1x_2x_3^2\geq
		4x_1x_2^3+4x_1x_3^3+4x_2x_3^3+4x_1^3x_3+12x_1^2x_2x_3+12x_1x_2^2x_3,
		$ \ \  or equivalently,$$(x_1+x_2-x_3)^4\geq8(x_1 x_2^3-x_2^3x_3),$$  with equality if and only if   $x_2=0$ and $x_1=x_3$;
		\item [(f)] $x_1^4+x_2^4+x_3^4+6x_1^2x_2^2+6x_1^2x_3^2+6x_2^2x_3^2+4x_1x_3^3+4x_2^3x_3+12x_1^2x_2x_3\geq4x_1^3x_2+4x_2x_3^3
		+4x_1x_2^3+4x_1^3x_3+12x_1x_2x_3^2+12x_1x_2^2x_3,
		$ \ \  or equivalently,$$(x_3+x_2-x_1)^4\geq8(x_2x_3^3-x_1 x_3^3),$$ with equality if and only if   $x_2=0$ and $x_1=x_3$;
		\item [(g)] $x_1^4+x_2^4+x_3^4+6x_1^2x_2^2+6x_1^2x_3^2+6x_2^2x_3^2+4x_2x_3^3+4x_1^3x_3+12x_1x_2^2x_3\geq
		4x_1x_2^3+4x_1x_3^3+4x_1^3x_2+4x_2^3x_3+12x_1^2x_2x_3+12x_1x_2x_3^2,
		$ \ \  or equivalently,$$(x_1+x_3-x_2)^4\geq8(x_1 x_3^3-x_2x_3^3),$$ with equality if and only if   $x_3=0$ and $x_1=x_2$.	\end{itemize} 
\end{theorem}

\section{Conclusions}

For a 4th-order 3-dimensional symmetric tensor with its entries  $1$ or $-1$,   the analytic  necessary and sufficient conditions are established for strict copositivity and copositivity,  respectively.  These conditions can be applied to  verify (strict) copositivity of a general 4th order 3-dimensional symmetric tensor. Several  (strict) inequalities  of ternary quartic homogeneous polynomial are built by means of these analytic conditions.

%\bibliographystyle{amsplain}

	
	
	
	% \end{thebibliography}

\end{document}